%% file: BoundedMonomialOperators.tex
\documentclass[12pt]{article}
\usepackage{verbatim}
\usepackage{amssymb}
\usepackage{amsthm}
\usepackage{amsmath}
\usepackage{mathrsfs}
\usepackage[margin=1in]{geometry}

\input def_latex3

\usepackage[usenames,dvipsnames] {color}

\definecolor{dark_purple}{rgb}{0.4, 0.0, 0.4}
\definecolor{dark_green}{rgb}{0.0, 0.7, 0.0}

\numberwithin{equation}{section}

\title{Monomial Operators}
\author{Jim Agler
and
John E. M\raise.5ex\hbox{c}Carthy
\thanks{Partially supported by National Science Foundation Grant  
DMS 2054199}
}
\date{\today}

\newcommand{\mc}{M\raise.45ex\hbox{c}Carthy}
\DeclareMathOperator{\re}{Re}

\newcommand\g{\gamma}

\usepackage{mathtools}

\newcommand\NN{{\mathbb N}}

\newcommand\hh{{\mathbb S}}
\newcommand\hhr{{\mathbb H}_{\rho}}
\newcommand\wt{\widehat{T}}
\newcommand\fc{f^{\cup}}
\newcommand\gc{g^{\cup}}
\newcommand\hc{h_u^{\cup}}
\newcommand\bc{\beta^{\cup}}
\newcommand\har{{\mathcal H}^2({\hh})}

\begin{document}

\bibliographystyle{plain}

\maketitle

\begin{abstract}
We study monomial operators on $ L^2[0,1]$, that is bounded linear operators that map each monomial $x^n$ to a multiple of $x^{p_n}$ for some $p_n$. We show that they are all unitarily equivalent to weighted composition operators on a Hardy space.
We characterize what sequences $p_n$ can arise. In the case that $p_n$ is a fixed translation of $n$, we give a criterion for boundedness of the operator.
\end{abstract}

\section{Introduction}

When studying polynomial approximation in $ L^2[0,1]$, the following class of operators arises naturally.
\begin{definition}
{\rm
A {\em monomial operator} is a bounded linear operator $T: L^2[0,1] \to L^2 [0,1]$ with the property  that there exist constants $c_n$ and $p_n$ so that
\be
\label{eqa1}
T : x^n \mapsto c_n \, x^{p_n} \qquad \forall n \in \NN.
\ee
If, in addition, there is some $\tau$ so that $p_n = n + \tau$ for every $n$, we call it
a {\em flat monomial operator}. }
\end{definition}
The powers $p_n$ may be complex, but must lie in the half plane 
\[
\hh := \{ s \in \C : \re(s) > -\frac{1}{2} \},
\]
in order for $x^{p_n}$ to lie in $L^2$.
Well-known examples of monomial operators include the
Hardy  operator $H$
\beq
 H : f &\mapsto& \frac{1}{x} \int_0^x f(t) dt \\
 H x^n & = & \frac{1}{n+1} x^n ;
\eeq
the operator  $M_x$ of multiplication by $x$; and the Volterra operator
\beq
V = M_x H :  f &\mapsto&  \int_0^x f(t) dt \\
 V x^n & = & \frac{1}{n+1} x^{n+1} .
\eeq
In \cite{amHI} the authors studied flat monomial operators, and showed that they  all
leave invariant every subspace of the form 
 \[
\{ f \in L^2[0,1] : f = 0 {\rm \ on\ }
[0,t] \} .
\]
  It was shown independently  by Brodskii \cite{br57} and Donoghue \cite{don57} that these  subspaces 
 are exactly the invariant subspaces for the Volterra operator.
  
 It is the purpose of this note to examine general monomial operators of the form \eqref{eqa1}.
To describe them, it is convenient to introduce a Hardy space associated with $\hh$.
The Hardy space of the unit disk, which we denote $H^2(\D)$, is the Hilbert space of holomorphic functions on 
$\D$ with finite norm, where the norm is given by
\[
\| \phi \|^2_{H^2(\D)} \= \sup_{0 < r < 1} \frac{1}{2\pi} \int_0^{2\pi} |\phi (r e^{i\theta})|^2 d \theta .
\]
 There are two distinct definitions for the Hardy space of a half-plane. Let $\lambda$  be the linear fractional transformation
give by
\be
\label{eqa5.4}
\lambda(s)  = \frac{s}{s+1},
\ee
 that takes $\hh$ onto $\D$. 
 \begin{definition}
 By $H^2(\hh)$ we mean $\{ \phi \circ \lambda : \phi \in H^2(\D) \}$, and the norm is defined so
 that $ \phi \mapsto \phi \circ \lambda$ is unitary.
\end{definition}
An equivalent definition of $H^2(\hh)$ is the set of functions $f$ that are holomorphic in $\hh$ and such that
$|f|^2$ has a harmonic majorant there. The norm squared is equal to the value of this harmonic majorant at $0$.
If $f \in H^2(\hh)$, then it has boundary values a.e. (see e.g. \cite{dur70}), and the norm is given by
\be
\label{eqa5.5}
\| f \|^2 \= \frac{1}{2\pi} \int_{-\i}^\i \frac{|f(-\frac{1}{2} + it)|^2}{t^2 + \frac{1}{4}} dt .
\ee
 The reproducing kernel for $H^2(\hh)$ is given by 
\be
\label{eqa6}
k(s,u) \= \la k_u, k_s \ra \= \frac{(1 + s)(1 + \bar u)}{1 + s + \bar u} .
\ee
In light of \eqref{eqa6}, the map 
\beq
U:  L^2[0,1] & \to & H^2(\hh) \\
 x^s & \mapsto & \frac{1}{1+s} k_{\bar s} 
\eeq
extends to a unitary from $L^2[0,1]$ onto $H^2(\hh)$. 
In \cite{amMS} we show that for any $f \in L^2$, we have
\[
Uf (s) = (1+s) \int_0^1 f(x) x^{s}dx .
\]

If $T$ is a monomial operator, we shall
define
\be
\label{eqa9}
\wt \ = \ U T U^* :  H^2(\hh) \to  H^2(\hh) .
\ee

One might wonder whether the definition of  monomial operator should  require that it take $x^s$ to some multiple of a monomial for every $s \in \hh$; our first theorem asserts that this always happens, just by assuming it on the natural numbers.
It also shows that, after moving to the Hardy space as in \eqref{eqa9}, monomial operators correspond to the adjoints of weighted composition operators.

If $g$ is a  holomorphic function on some domain, let us define
\[
\gc(s) \= \overline{g(\overline{s})} .
\]
By \eqref{eqa5.5} we see that if $f \in H^2(\hh)$, then  $\| \fc \|_{H^2(\hh)} = \| f \|_{H^2(\hh)}$.
\bt
\label{thma1}
Let $T:  L^2[0,1] \to L^2 [0,1]$ be a monomial operator given by \eqref{eqa1}.
Then there exists a holomorphic map $\beta: \hh \to \hh$ and a function $g \in H^2(\hh)$ so
that, for every $s \in \hh$, we have
\be
\label{eqa7}
T (x^s) \=  \frac{1+ \beta(s)}{1 + s} g(s)  \ x^{\beta(s)}.
\ee
Moreover, we have
\be
\label{eqa8}
\wt^* f (s) \= \gc(s) f (\bc (s)) \qquad \forall f \in H^2(\hh).
\ee
\et
Equation \eqref{eqa8} says
\[
\wt^* \= M_{\gc} C_{\bc} ,
\]
where $C_\beta$ denotes the composition operator $f \mapsto f \circ \beta$, and $M_{g}$  denotes the multiplication
operator $ f \mapsto g f$.

Weighted composition operators have been studied for some time. See
e.g. \cite{fo64, coga06, krmo07, coko10, bona10, acks18, cp21}.
By Littlewood's subordination principle, the operator $C_\beta$ is always bounded whenever $\beta$ is a holomorphic self-map of $\hh$ \cite[Thm. 10.4.2]{zhu}. A multiplication operator $M_h$ is bounded if and only if $h \in H^\i(\hh)$.
It was observed in \cite{guna08} that 
it is possible
for the product $M_h C_\beta$ to be bounded even when $M_h$ is not.

A consequence of Theorem \ref{thma1}
is that it allows us to specify what sequences $p_n$ can occur in \eqref{eqa1}.
\begin{corollary}
\label{cora1}
Let $ (p_n)$ be a sequence in $\hh$. Then there exists some choice of scalars $c_n$, not all zero, so that
$T : x^n \mapsto c_n x^{p_n}$ extends to be a bounded linear operator on $L^2[0,1]$ if and only if 
\[
\left[ \frac{p_m + \overline{p_n} + 1}{ m + n + 1} \right] \ \geq \ 0 .
\]
\end{corollary}

In Section \ref{secflat} we study which pairs of functions $\beta$ and $g$ 
give rise to a bounded operator in \eqref{eqa7}.
We answer this question only for flat monomial operators, i.e. when $\beta(s) = s + \tau$ for some constant $\tau$. If $\Re (\tau) < 0$, it follows from Corollary \ref{cora1} that such a $T$ can never be bounded.

If $\Re(\tau) \geq 0$,  we have: 
\bt
\label{thma3}
Let $\Re (\tau) \geq 0$, and let $T$ be defined  by
\[
T (x^n) \=  \frac{1+ n + \tau}{1 + n} g(n)  \ x^{n + \tau}
\]
for some function $g \in H^2(\hh)$.

(i) If $\Re(\tau) > 0$, then
 $T$ extends to be a bounded linear operator from $L^2[0,1]$ to $L^2[0,1]$
if and only if
 the Poisson integral of $|g|^2$ is bounded on
all half-planes that are strictly contained in $\hh$.

(ii) If $\Re(\tau) = 0$, then $T$ is bounded if and only if $g$ is bounded on $\hh$.

\et
 This theorem is proved in Theorems \ref{thmc1} and  \ref{thmc2}.

\section{Proof of Theorem \ref{thma1}}
\label{secb}

{\sc Proof of Theorem \ref{thma1}.}
Step 1. Let $\wt$ be given by \eqref{eqa9}.
Define $g$ by $\gc = \wt^* 1 $.
We have
\beq
\wt k_{ n} &\=& U T (1+n) x^n \\
&=& (1+n) U c_n x^{p_n} \\
&=& c_n \frac{1+n}{1+ p_n} k_{\overline{p_n}} .
\eeq
Since $k_0 = 1$, we have
\beq
\la \gc , k_n \ra
&\=&
\la 1,  c_n \frac{1+n}{1+ p_n} k_{\overline{p_n}} \ra\\
&=& \overline{c_n}  \frac{1+n}{1+\overline{ p_n}}.
\eeq
This gives
\[
g(n) \=  c_n \frac{1+n}{1+ p_n} , 
\]
and so we can write
\be
\label{eqb1}
\wt k_n \= g(n) k_{\overline{p_n}}.
\ee
For any $u \in \hh$, let $h_u = \wt^* k_{ u}$.
We get from \eqref{eqb1}
\begin{eqnarray}
\nonumber
\la h_u, k_n \ra &\=& 
\gc(n) \la k_{ u}, k_{\overline{p_n}}\ra\\
\label{eqb1.5}
\att
\Rightarrow\quad h_u(n) &=& \gc(n) \frac{( 1+  \bar u)( 1 + \overline{p_n})}{1 + \bar  u + \overline{p_n}}\\
\label{eqb2}
\Rightarrow\quad \hc(n) &\= &
g(n) \frac{( 1+  u)( 1 + {p_n})}{1 +  u + {p_n}}.
\end{eqnarray}
\att
Define $\beta$ by
\be
\label{eqb3}
\beta(s) \= \frac{(1+ u) (\hc(s) - g(s))}{(1+ u) g(s) - \hc(s)} .
\ee
As $g$ and $h_u$ are both in the Hardy space $H^2(\hh)$, a priori we know that $\beta$ is in the
Nevanlinna class of meromorphic functions on $\hh$, the class of quotients of $H^2$ functions.
Moreover, it follows from \eqref{eqb2} that 
\[
p_n \=  \frac{(1+ u) (\hc(n) - g(n))}{(1+ u) g(n) - \hc(n)}  \= \beta(n) .
\]
Observe that $\NN$ is not a zero set for $H^2(\hh)$---indeed, $(\lambda(n) = \frac{n}{n+1})$ is not a Blaschke sequence for $H^2(\D)$. Therefore it is a set of uniqueness for the Nevanlinna class, and hence 
$\beta$ is the unique function in the class that satisfies $\beta(n) = p_n$.
In particular, the definition \eqref{eqb3} is actually independent of $u$.

We can write \eqref{eqb1.5} as
\be
\label{eqb4}
(\wt^* k_u )(n) \= \gc(n) k_u(\bc(n)) \qquad \forall n \in \NN.
\ee
Since $\NN$ is a set of uniqueness, \eqref{eqb4} holds everywhere
\be
\label{eqb5}
(\wt^* k_u )(s) \= \gc(s) k_u(\bc(s)) \qquad \forall s \in \hh.
\ee
Step 2. We must show that $\beta: \hh \to \hh$.
Suppose there is some point $s \in \hh$ where $g(s) \neq 0$ and $\bc(s) = w$ is in $\C \setminus \hh$.
Then there is a sequence $q_n$ such that each $q_n$ is a finite linear combination of kernel functions,
$\| q_n \| \leq 1$, and $q_n(w) \to \i$. By \eqref{eqb5}, we have
\be
\label{eqb6}
(\wt^* q_n) (s) \= \gc(s) q_n(w) .
\ee
The right-hand side of \eqref{eqb6} tends to infinity, the left-hand side is bounded by $\| T \| \| k_s \|$, a contradiction.
Since $\bc$ is meromorphic, we conclude that whenever $\gc(s) \neq 0$ then $\beta(s) \in \hh$.
Therefore any singularities of $\bc$ on the zero set of $\gc$ are removable, so we conclude that $\bc$,
and hence also $\beta$, is a self-map of $\hh$.

Thus we have proved \eqref{eqa8} for any finite linear combination of kernel functions, and so, 
by a limiting argument, it is true for all $f \in H^2(\hh)$.

Step 3. From \eqref{eqb5}, we have
\beq
\la \wt k_s , k_u \ra &=&  \overline{\la \wt^* k_u, k_s \ra}  \\
 &\=& g(\bar s) \overline{k_u(\bc(s))} \\
 &=& g(\bar s) \frac{(1+u)(1 + \beta(\bar s))}{1 + u + \beta(\bar s))} \\
 &=& g(\bar s) k_{\bc(s)} (u).
\eeq
So 
\be
\label{eqb7}
\wt k_s \= g(\bar s) k_{\bc(s)} .
\ee
Therefore
\beq
T x^s &\=& U^* \wt U [x^s] \\
&=& U^* \wt [\frac{1}{1+s} k_{\bar s}] \\
&=& 
U^* [ \frac{1}{1+s} g(s) k_{\overline{\beta (s)}}] \\
&=& 
\frac{1+ \beta(s)}{1+s} g(s) x^{\beta (s)} .
\eeq
This proves \eqref{eqa7}.
\ep

{\sc Proof of Corollary \ref{cora1}.}
By Theorem \ref{thma1}, a necessary condition for the
existence of a non-zero bounded $T$ that maps each $x^n$ to a multiple of $x^{p_n}$ is that there be some
holomorphic self-map $\beta$ of $\hh$ that maps $n$ to $p_n$. This condition is also sufficient, since
choosing
\[
c_n \= \frac{1+p_n}{1+n}
\] 
gives $U T U^*$ is the adjoint of $C_{\bc}$, which is bounded.

When is there a map $\beta: \hh \to \hh$ that interpolates $n$ to $p_n$?
Composing with the Riemann map $\lambda$ from \eqref{eqa5.4}, this is equivalent to asking when
there exists $\phi = \lambda \beta \lambda^{-1}$ from $\D$ to $\D$ that maps $\frac{n}{n+1}$ to
$\frac{p_n}{p_n +1}$.
By Pick's theorem \cite[Thm. 1.81]{amy20}, this occurs if and only if
\be
\label{eqb8}
\left[ \frac{ 1 - \frac{p_m}{p_m +1}\frac{\overline{p_n}}{\overline{p_n} +1}}{1 - \frac{m}{m+1} \frac{ \bar n}{\bar n +1}} \right] \ \geq \ 0 .
\ee
Rearranging \eqref{eqb8}, we get
\[
\left[ \frac{m+1}{p_m + 1} \frac{ \bar n +1}{\overline{p_n} + 1} 
\frac{p_m + \overline{p_n} + 1}{m + \bar n + 1} \right] \ \geq \ 0.
\]
As conjugating by the rank one operator 
\[
\left[ \frac{m+1}{p_m + 1} \frac{ \bar n +1}{\overline{p_n} + 1} \right]
\]
does not affect positivity, and $n = \bar n$, we get that an interpolating $\beta$ exists if and only if
\[
\left[ \frac{p_m + \overline{p_n} + 1}{m +  n + 1} \right] \ \geq \ 0.
\]
\ep

\section{Flat Monomial Operators}
\label{secflat}

In this section, we shall consider operators of the form 
\be
\label{eqc2}
T : x^n \mapsto c_n x^{n + \tau} ,
\ee
and in particular when they are bounded.
Note first that if $\tau = \tau_0 + i \tau_1$, with $\tau_0$ and $\tau_1$ real, then the effect of $\tau_1$ is to multiply everything in the range by $x^{i\tau_1}$, which is unimodular. So without loss of generality we can assume that $\tau$ is real. Moreover, by Corollary \ref{cora1}, $T$ can only be bounded if $\tau \geq 0$.

Let us first handle the case $\tau = 0$. In the terminology of  Theorem \ref{thma1}, $\beta(s) =s$, so  $T$ is bounded if and only if there is a bounded function $g$ on $\hh$ that satisfies $g(n) = c_n$.
By Pick's theorem, this happens if and only if
\be
\label{eqc1}
\left[ \frac{1-c_m \overline{c_n}}{1 + m + n} \right] \ \geq \ 0 .
\ee
So we get:
\bt
\label{thmc1}
The map $T: x^n \mapsto c_n x^n$ extends to be a bounded linear map if and only if \eqref{eqc1} holds.
Moreover, if it is non-zero, it is never compact.
\et

When $\tau > 0$, how do we determine whether \eqref{eqc2} extends to be bounded? 
Let
\be
\label{eqc2.1}
\gamma(s) \= s + \tau .
\ee
 We must find $g \in H^2 (\hh)$ that satisfies 
\[
g(n) \= \overline{c_n} \ \frac{1+n}{1 + n + \tau} .
\]
Then $T$ is bounded if and only if $M_{g} C_\gamma$ is bounded.
(We have changed notation slightly from Theorem \ref{thma1} to avoid the use of $\gc$ and
to emphasize that we have a fixed choice of $\beta$).
To investigate $M_g C_\gamma$, we turn to the Poisson kernel.

The Poisson kernel for $\hh$ at a point $s = -\frac{1}{2} + \sigma + i t$ is given by
\[
P_{ -\frac{1}{2} + \sigma + i t} (-\frac{1}{2} + iy) \= \frac{1}{\pi} \frac{\sigma}{\sigma^2 + (y-t)^2} .
\]
The Poisson integral of some function $f$ defined on the line $\{ \Re (s) = - \frac{1}{2} \}$ is given by
\[
P[f] (-\frac{1}{2} + \sigma + it)
\=
\frac{1}{\pi} \int_{-\i}^\i \frac{\sigma}{\sigma^2 + (y-t)^2} f(-\frac{1}{2} + iy) dy
\]
It is convenient to introduce another Hardy space, $\har$.
\begin{definition}
The space $\har$ consists of all functions $f$ that are holomorphic in $\hh$ and satisfy
\[
\sup_{ x > - \frac{1}{2}} \int_{-\i}^\i |f(x+iy)|^2 dy < \i .
\]
\end{definition}
Every function $f \in \har$ has non-tangential boundary limits almost everywhere, and its norm is then
given by (see \cite{gar81}):
\[
\| f \|^2_{\har} \=  \int_{-\i}^\i |f(-\frac{1}{2}+iy)|^2 dy .
\]
Let us write
\[
\hhr := \{ s \in \C : \Re(s) > \rho \} .
\]
So $\hh = {\mathbb H}_{\frac{1}{2}}$ in this notation.

\bt
\label{thmc2}
The operator $M_{g} C_\gamma$ is bounded if and only if the Poisson integral of $|g|^2$
is bounded on some (and hence every) half-plane $\hhr$ for $ \rho > - \frac{1}{2}$.
\et

{\sc Proof of Thm. \ref{thmc2}.}
The map
\[
W: f(s) \mapsto \frac{1}{\sqrt{2\pi}} \frac{1}{s+1} f(s) 
\]
is a unitary from $H^2(\hh)$ onto $\har$. One checks 
\be
\label{eqc3}
W M_g C_\gamma W^* \= M_{\frac{s+\tau + 1}{s+1} g(s)} C_\gamma .
\ee
For $s$ in $\hh$ we have
\[
1 \leq | \frac{s+\tau + 1}{s+1} | \leq 1 + 2\tau,
\]
so for boundedness purposes we can drop this factor in \eqref{eqc3} and conclude
that $M_g C_\gamma$ is bounded on  $H^2(\hh)$ if and only if it is bounded on $\har$.
Thus we wish to know for which $g$ does
\[
\int_{-\i}^\i |g(-\frac{1}{2}+iy)|^2  |f(-\frac{1}{2}+ \tau +iy)|^2   dy 
\ \lesssim \ 
\int_{-\i}^\i |f(-\frac{1}{2}+iy)|^2 dy 
\] 
hold?
This is the same as asking when integrating along the vertical line $\{ \Re(s) = -\frac{1}{2} + \tau \}$
with weight $|g(-\frac{1}{2}+iy)|^2$ is a Carleson measure for $\har$.
By \cite[VI.3]{gar81}, this happens if and only if
\be
\label{eqc4}
\sup_{\sigma > 0, t \in \R} \int_{-\i}^\i \frac{\sigma}{(\tau + \sigma)^2 + (y-t)^2} 
|g(-\frac{1}{2}+iy)|^2 dy \ = C \ < \ \i .
\ee
Taking $\sigma = \tau$ in \eqref{eqc4}, we get that $M_g C_\gamma$ is bounded implies 
\[
\sup_{ t \in \R} \int_{-\i}^\i \frac{\tau}{(2\tau )^2 + (y-t)^2} 
|g(-\frac{1}{2}+iy)|^2 dy \ \leq \ K ,
\]
for some constant $K$,
 so
\[
\sup_{ t \in \R} \int_{-\i}^\i \frac{\tau}{(\tau )^2 + (y-t)^2} 
|g(-\frac{1}{2}+iy)|^2 dy \ \leq \ 4K .
\]
If $u= -\frac{1}{2} + \tau + \sigma + i t$ is any point in the half-plane ${\mathbb H}_{-\frac{1}{2} + \tau}$,
we have
\[
P_u(-\frac{1}{2} + iy) \= \frac{1}{\pi} \frac{\sigma + \tau}{(\sigma + \tau)^2 + (y-t)^2} .
\]
As
\[
\frac{\sigma + \tau}{(\tau + \sigma)^2 + (y-t)^2} \ \leq 
\frac{\sigma }{(\tau + \sigma)^2 + (y-t)^2}
+ \frac{ \tau}{(\tau )^2 + (y-t)^2},
\]
we conclude that \eqref{eqc4} implies that the Poisson integral of $|g|^2$ is bounded on
${\mathbb H}_{-\frac{1}{2} + \tau}$.
Conversely, the boundedness of 
the Poisson integral of $|g|^2$ on ${\mathbb H}_{-\frac{1}{2} + \tau}$ implies \eqref{eqc4}.

Finally, to see that boundedness of the Poisson integral of $|g|^2$ on ${\mathbb H}_{-\frac{1}{2} + \sigma}$
implies boundedness on ${\mathbb H}_{-\frac{1}{2} + \rho}$ for any $ 0 < \rho < \sigma$,
we just observe that
\[
\frac{1}{\pi} \frac{\rho}{\rho^2 + (y-t)^2}
\ \leq \  ( \frac{\sigma^2}{\rho^2} )
\frac{1}{\pi} \frac{\sigma}{\sigma^2 + (y-t)^2}.
\]
\ep

\begin{exam}
Let $\tau > 0$, and take $g(s) = \frac{1}{(1+s) (s+\frac{1}{2})^c}$ for some $ 0 < c < \frac{1}{2}$.
Then $g$ is not bounded on $\hh$, so $M_g$ is not bounded.
But
\beq
P[|g|^2] (-\frac{1}{2} +\sigma + it) 
&\=&
\frac{1}{\pi} \int_{-\i}^{\i}  \frac{\sigma}{\sigma^2 + (y-t)^2} \frac{1}{|y|^{2c}} \frac{1}{\frac{1}{4} + y^2}dy \\
&\leq& \frac{1}{\pi \sigma} \int_{-\i}^{\i}   \frac{1}{|y|^{2c}} \frac{1}{\frac{1}{4} + y^2} dy
\\
&\leq& \frac{1}{\pi \sigma} \left[ \int_{-1}^1 \frac{4}{|y|^{2c}} dy + \int_{-\i}^\i  \frac{1}{\frac{1}{4} + y^2}  dy \right] \\
&\=&
 \frac{1}{\pi \sigma} \left[ \frac{8}{1-2c}  + 2\pi \right]
\eeq
This is bounded in each half-plane ${\mathbb H}_\rho$, so by Theorem \ref{thmc2} $M_g C_\gamma$ is bounded.
\end{exam}

\section{Unitary Monomial Operators}

Bourdon and Narayan characterized unitary weighted composition operators \cite{bona10}.
Their theorem (translated to $\hh$) is:
\bt
\label{thmd2} (Bourdon-Narayan)
The operator $M_g C_\beta$ is unitary on $H^2(\hh)$ if and only if 
$\beta$ is an automorphism of $\hh$ and $g(s) =  e^{i\theta} k_{s_0}/\| k_{s_0} \|$,
where $s_0 = \beta^{-1} (0)$.
\et
We shall give a  proof of  their theorem directly in the context of monomial operators.
First, let us describe the automorphisms of $\hh$ in a convenient way.
\begin{lemma}
\label{lemd1}
The function $\beta: \hh \to \C$ is an automorphism of $\hh$ if and only if there is a function $\phi : \hh \to \C$ so that
\be
\label{eqd2}
\frac{1 + \beta(s) + \overline{\beta(t)}}{1 + s + \overline{t}} \= \phi(s) \overline{\phi (t)} .
\ee
\end{lemma}
\bp
Notice that $\beta$ is an automorphism of $\hh$ if and only if $ b := \lambda \circ \beta \circ \lambda^{-1}$ is an automorphism of $\D$.
By Pick's theorem (see \cite[Sec. 2.6]{amy20}) the latter occurs if and only if
\be
\label{eqd1}
\frac{ 1 - \overline{b(w)} b(z)}{ 1 - \bar w z} \ = \ \overline{\psi(w)} \psi (z) 
\ee
for some function $\psi$ on $\D$.
Doing some algebra, \eqref{eqd1} becomes \eqref{eqd2} with
\[
\phi(s) \= \frac{1 + \beta(s)}{1+s} \psi( \frac{s}{1+s} ) .
\]
\ep
Consequently, we have the following characterization of unitary monomial operators (which can also be derived from Bourdon-Narayan and Theorem \ref{thma1}).
\bt
The operator 
\[
T : x^s \mapsto c(s) x^{\beta(s)} 
\]
is unitary on $L^2[0,1]$ if and only if $\beta$ is an automorphism of $\hh$ and $c(s)$ is defined by 
\be
\label{eqd3}
c(s) \= \frac{e^{i\theta}}{\sqrt{1 + 2 \Re \beta(0)}} \frac{ 1 + \overline{\beta(0)} + \beta(s)}{1 + s } .
\ee
\et 
\bp
Since $\beta$  is holomorphic by Theorem \ref{thma1} and 
has to be non-constant for $T$ to be bounded, the range of $T$ must be dense. Therefore it is unitary if and only if it is isometric.
It is isometric if and only if for every $s,t$ we have
\beq
\la x^s , x^t \ra &\=& \la c(s) x^{\beta(s)}, c(t) x^{\beta(t)} \ra \\
\Leftrightarrow  \quad \frac{1}{1+ s + \bar t} & \= & \frac{ c(s) \overline{c(t)} }{1 + \beta(s) +  \overline{\beta(t)}}.
\eeq
That means \eqref{eqd1} holds, so $\beta$ is an automorphism.  Letting $t=0$ 
we get
\[
c(0) \= e^{i \theta} \sqrt{1 + 2 \Re \beta(0)},
\]
and solving for $c (s)$ we get
\eqref{eqd3}.

Conversely, suppose \eqref{eqd3} holds and $\beta$ is an automorphism. By Lemma \ref{lemd1}, we have
\eqref{eqd2} for some $\phi$. Letting $t=0$ and solving, we get that $\phi$ is given by \eqref{eqd3}.
\ep

\section{Questions}

\begin{question}
If $M_g C_\g$ is bounded, with $\gamma$ as in \eqref{eqc2.1}, can one approximate it in norm by 
operators of the form  $M_{g_n} C_\g$ where $g_n \in H^\i(\hh)$?
\end{question}

\begin{question}
A generalization of the previous question is for non-flat monomial operators.  Can every bounded operator of
the form $M_g C_\beta $ be approximated by operators $M_{g_n} C_\beta$ with bounded $g_n$?
\end{question}

\begin{question} 
What are necessary and sufficient conditions for the functions $\beta$ and $g$ so that the operator $T$ defined by 
\eqref{eqa7} is bounded?
Compact?
\end{question}

%
%
%

On behalf of all authors, the corresponding author states that there is no conflict of interest.

\bibliography{references_uniform_partial}
\end{document}

%% file: def_latex3.tex
\usepackage{latexsym}
\usepackage{amssymb}
\usepackage{euscript}

\def\mcc{M\raise.5ex\hbox{c}C}
\def\mccarthy{M\raise.5ex\hbox{c}Carthy}

\let\i=\infty

\def\la{\langle}
\def\ra{\rangle}
\def\={\ = \ }

\def\C{\mathbb C}
\def\R{\mathbb R}

\def\D{\mathbb D}

\def\be{\setcounter{equation}{\value{theorem}} \begin{equation}}
\def\ee{\end{equation} \addtocounter{theorem}{1}}
\def\beq{\begin{eqnarray*}}
\def\eeq{\end{eqnarray*}}
 
\def\att{\addtocounter{theorem}{1}}

\def\bp{{\sc Proof: }}
\def\ep{{}{\hfill $\Box$} \vskip 5pt \par}

\def\bl{\begin{lemma}}
\def\el{\end{lemma}}
\def\bt{\begin{theorem}}
\def\et{\end{theorem}}
\def\bprop{\begin{prop}}
\def\eprop{\end{prop}}
\def\bd{\begin{definition}}
\def\ed{\end{definition}}
\def\br{\begin{remark}}
\def\er{\end{remark}}
\def\bexer{\begin{exercise}}
\def\eexer{\end{exercise}}

\newtheorem{theorem}{Theorem}[section]

\newtheorem{lemma}[theorem]{Lemma}
\newtheorem{corollary}[theorem]{Corollary}

\newtheorem{definition}[theorem]{Definition}

\newtheorem{question}[theorem]{Question}